\documentclass{amsart}[12pt]
\usepackage{amsmath}
\usepackage{amstext}
\usepackage{amssymb}
\usepackage{amsthm}
\usepackage{amscd}

\usepackage{lscape}
\usepackage{longtable}

\usepackage{epsfig}
\input txdtools
\let\et=\etexdraw
\def\etexdraw{\drawbb\et}

\def\baselinestretch{1.5}
\theoremstyle{plain}
\newtheorem{thm}{Theorem}[section]

\newtheorem{prop}[thm]{Proposition}

\theoremstyle{definition}

\theoremstyle{remark}

\DeclareMathOperator{\initial}{in}

\DeclareMathOperator{\HF}{HF}
\DeclareMathOperator{\HS}{HS}
\DeclareMathOperator{\HP}{HP}

\begin{document}

\def\Ker {\operatorname{Ker}\nolimits}
\def\Im {\operatorname{Im}\nolimits}
\def\Image {\operatorname{Image}\nolimits}
\def\Syz {\operatorname{Syz}\nolimits}
\def\initial {\operatorname{in}\nolimits}
\def\Gin{\operatorname{Gin}\nolimits}
\def\Spec{\operatorname{Spec}\nolimits}
\def\D{\operatorname{D}\nolimits}
\def\V{\operatorname{V}\nolimits}
\def\H {\operatorname{H}\nolimits}
\def\E {\operatorname{E}\nolimits}
\def\V {\operatorname{V}\nolimits}
\def\nE {\operatorname{e}\nolimits}
\def\nV {\operatorname{v}\nolimits}
\def\C {\operatorname{\cal C}\nolimits}

\newcommand {\Reduce}[3] {#1\xrightarrow[#2]{} #3}

\title
[Counting monomials]
{Counting monomials}
\author{Mordechai Katzman}
\address{Department of Pure Mathematics,
University of Sheffield, Hicks Building, Sheffield S3 7RH, United Kingdom\\
{\it Fax number}: 0044-114-222-3769}
\email{M.Katzman@sheffield.ac.uk}




\begin{abstract}
This paper presents two enumeration techniques based on Hilbert
functions. The paper illustrates these techniques by solving two
chessboard problems.
\end{abstract}

\maketitle

\section{Introduction and preliminaries.}

The purpose of this note is to illustrate two powerful enumeration techniques
based on computational Commutative Algebra methods.

By way of illustration I chose to apply these methods to the following two elementary problems:
\begin{enumerate}
\item
Consider a $n\times n$ chessboard.
What is the maximal number of unattacked squares in the board after placing on it $k$ queens?
More generally, in how many ways can we place $k$ queens on a chess board to obtain exactly $u$ unattacked squares?

\item
Consider an infinite chessboard.
How many squares can a knight reach in $d$ moves?
How many squares can be reached in $d$ moves and no less?
\end{enumerate}
Although these problems are phrased in the language of chess,
they are specific instances of more general graph-theoretical problems.
The enumeration techniques presented here answer these more general problems.

At the heart of the methods presented in this paper
are the notions of graded modules and their Hilbert functions.
In essence, we will reduce each of the problems above to a problem about the enumeration of
sets of monomials, and this enumeration will be achieved using Hilbert functions.

While the application of Hilbert functions to the problems presented in this paper is new,
the use of Hilbert functions in combinatorics is not.
The solution of some simple enumeration problems using Hilbert functions,
such as finding the independence number of a graph,
has long been part of the folklore of computational commutative algebra experts.
An early and striking example of the use Hilbert functions in combinatorics
is Richard P.~Stanley's work on magic squares
(I refer the reader to \cite{R} for an accessible and thoroughly enjoyable account of this work.)

We now review graded modules and Hilbert functions.
Throughout this paper, all rings are commutative and with $1$; $K$ will always denote a field.

A $K$-algebra $R$ is $\mathbb{N}^N$-graded if we can write
$$R=\bigoplus_{\mathbf{a}\in \mathbb{N}^N} {R}_\mathbf{a} ,$$
a direct sum of abelian groups, and the direct summands satisfy
$${R}_\mathbf{a} {R}_\mathbf{b} \subseteq {R}_{\mathbf{a}+\mathbf{b}}$$
for all $\mathbf{a},\mathbf{b} \in \mathbb{N}^N$.
Henceforth we shall also impose the condition
$R_\mathbf{0}=K$, which implies that each ${R}_\mathbf{a}$ is a $K$-vector space and that,
if $R$ is a finitely generated $K$-algebra,
each ${R}_\mathbf{a}$ is a finite dimensional $K$-vector space.
For each $\mathbf{a}\in \mathbb{N}^N$ we shall refer to the elements of
${R}_\mathbf{a}$ as being homogeneous of degree $\mathbf{a}$.

A fundamental example of such a graded $K$-algebra is the ring of polynomials $R=K[x_1, \dots, x_n]$.
We can endow $R$ with different graded structures.
We are all familiar with the $\mathbb{N}$-grading
$$R=\bigoplus_{a\in \mathbb{N}} {R}_a $$
in which each $R_a$ consists of the homogeneous polynomials of degree $a$.
We can define another grading as follows:
let $\mathbf{d}_1, \dots, \mathbf{d}_n\in \mathbb{N}^N$ and define the degree of a
monomial $x_1^{\alpha_1} \dots x_n^{\alpha_n}$ to be ${\alpha_1}\mathbf{d}_1+ \dots {\alpha_n}\mathbf{d}_n$.
We can now write
$$R=\bigoplus_{\mathbf{a}\in \mathbb{N}^N} {R}_\mathbf{a} ,$$
where each $\mathbb{R}_\mathbf{a}$ is the $K$-vector space spanned by all monomials of
degree $\mathbf{a}\in \mathbb{N}^N$.

Let $R$ be a $\mathbb{N}^N$-graded $K$-algebra. An $R$-module $M$ is graded if it has a $\mathbb{N}^N$-grading compatible with that of $R$, i.e.,
if we can write
$$M=\bigoplus_{\mathbf{a}\in \mathbb{N}^N} {M}_\mathbf{a} ,$$
a direct sum of abelian groups, and the direct summands satisfy
$${R}_\mathbf{a} {M}_\mathbf{b} \subseteq {M}_{\mathbf{a}+\mathbf{b}}$$
for all $\mathbf{a},\mathbf{b} \in \mathbb{N}^N$.

If $R$ is a polynomial ring as in the examples above and $I\subset R$ is a homogeneous ideal,
i.e., an ideal generated by homogeneous elements, then $R/I$ has a natural structure of a graded $R$-module.

\bigskip
Let $R$ be a $\mathbb{N}^N$-graded $K$-algebra and let $M$ be a graded $R$-module.
We define the Hilbert function $\HF_M$ of $M$ to be the function
$\HF_M : \mathbb{N}^N \rightarrow \mathbb{N}$ defined by
$\HF_M (\mathbf{a}) = \dim_K M_\mathbf{a}$.
The Hilbert series $\HS_M(t_1, \dots, t_N)$ of $M$ is the generating function of the Hilbert function, i.e.,
$$\HS_M(t_1, \dots, t_N) =\sum_{\mathbf{a}\in \mathbb{N}^N} \HF_M(\mathbf{a}) t_1^{a_1} \dots t_N^{a_N} .$$

If $R$ is a polynomial ring as in the examples above with its familiar $\mathbb{N}$-grading, and if we view $R$ as a graded $R$-module,
then $\HF_R (a)$ is just the number of monomials of degree $a$ in $n$ variables, i.e., $\HF_R (a)=\binom{a+n-1}{a}$, and
$\HS_R(t)=1/(1-t)^n$. If we were to assign degrees $\mathbf{d}_1, \dots, \mathbf{d}_n\in \mathbb{N}^N$ to $x_1, \dots, x_n$ we would
obtain
$$\HS_R(t_1, \dots, t_N)=\frac{1}{\prod_{i=1}^n 1-t_1^{\mathbf{d}_{i 1}} \dots t_N^{\mathbf{d}_{i N}}} .$$

Take $R$ to be a polynomial ring with its familiar $\mathbb{N}$-grading, let $I\subset R$ be a homogeneous ideal
and write $S=R/I$.
One can show that $\HF_S (a)$ is of polynomial type, i.e., it agrees with a polynomial, the Hilbert polynomial $\HP_S(a)$ of $S$, for all $a\gg 0$.
The degree of $\HP_S$ is one less than the Krull dimension of $S$. Also, one can write
$$\HS_S(t)=\frac{P(t)}{(1-t)^d}$$
where $P(t)$ is a polynomial which
does not vanish
at $t=1$ and $d$ is the Krull dimension of $S$.

\section{Unattacked squares}

We now consider the first question mentioned in the introduction. We naturally identify the squares of the $n\times n$ chessboard
with pairs $(i,j)$ where $1\leq i,j \leq n$.

We fix $n$, the size of the board.
Let $K$ be any field and define $R$ to be the polynomial ring in $2n^2$ variables
$$R=K[x_{11}, \dots, x_{nn},y_{11},\dots,y_{nn}] .$$
We assign degree $(1,0)$ to all the $x$ variables and degree $(0,1)$ to all the $y$ variables.

Roughly, the $x$ variables will correspond to squares in our $n\times n$ chessboard
which are occupied by queens while the $y$ variables will correspond to
unattacked squares on the board.

We define $I$ to be the ideal of $R$ generated by the squares of all variables together with
$$  \left\{ x_{ij} y_{lm} \,|\, \mathrm{a\ queen\ can\ move\ from\ square\ }(i,j)\mathrm{\ to\ square\ }(l,m) \right\} .$$
Notice that $I$, as any other ideal generated by monomials, is homogeneous with respect to the $\mathbb{N}^2$-grading of $R$.

For any $k>0$ define
$$\mu(k)=\max\{ \mu\in \mathbb{N} \, |\, \dim_K \left(R/I\right)_{(k,\mu)}>0 \} .$$

\begin{prop}
$\mu(k)$ is the maximal number of squares on the $n\times n$ chessboard which can remain unattacked after placing on it $k$ queens.
\end{prop}
\begin{proof}
Consider any monomial $M=\mathbf{x}^\alpha \mathbf{y}^\beta$ in $R$ whose image in $R/I$ is not zero.
Since $I$ contains the squares of all the variables, $M$ must be square-free and we may write
$$ M=x_{i_1,j_1} \cdot \ldots \cdot x_{i_\lambda,j_\lambda} y_{l_1,m_1} \cdot \ldots \cdot y_{l_\nu,m_\nu} .$$
where all the variables in this expression are distinct.
We next observe that for any $1\leq \xi \leq \lambda$ and $1\leq \zeta \leq \nu$,
a queen cannot move from square $(i_\xi, j_\xi)$ to square $(l_\zeta, m_\zeta)$, otherwise,
$x_{i_\xi, j_\xi} y_{l_\zeta, m_\zeta}$ would be one of the generators of $I$ and $M$ would be zero modulo $I$.
We showed that every monomial of degree $(\lambda, \mu)$ whose image in $R/I$ is not zero corresponds
to a configuration on the chessboard where the squares
$(i_1, j_1), \dots, (i_\lambda, j_\lambda)$
are occupied by queens and the squares
$(l_1, m_1), \dots, (l_\nu, m_\nu)$ are not attacked by any of these queens.

It is easy to see that the converse is also true and so we have established a bijection between
the configurations of $\lambda$ queens and $\nu$ unattacked squares and the set of
monomials of degree $(\lambda, \nu)$ which are not zero modulo $I$.

Notice that all the graded components
$\left(R/I\right)_{(\lambda, \nu)}$
are spanned as $K$-vector spaces by monomials of degree $(\lambda, \nu)$, and that a basis for
$\left(R/I\right)_{(\lambda, \nu)}$ is given by the set of all such monomials whose images in $R/I$ are not zero.
So now we can see that the condition
$$\dim_K \left(R/I\right)_{(k, \mu)}>0, \quad \dim_K \left(R/I\right)_{(k, \mu+1)}=0$$
can be translated using the bijection established above to the statement that
it is possible to place $k$ queens on the chessboard so that one can find $\mu$ unattacked squares but not $\mu+1$ unattacked squares.
\end{proof}

We now address
the more general question:
in how many ways $\Phi(k,u)$ can we place $k$ queens on a chessboard to obtain exactly $u$ unattacked squares?

\begin{prop}
For any $0\leq u \leq \mu(k)$
$$\Phi(k,u)=\HF_{R/I} (k, u) - \sum_{v=u+1}^{\mu(k)} \binom{v}{u} \Phi(k,v) .$$
\end{prop}
\begin{proof}
We proceed to prove this by reverse induction of $u$. When $u=\mu(k)$
the equality $\Phi(k, \mu(u))=\HF_{R/I} (k, \mu(u))$
follows easily from the discussion in the proof of the previous proposition.

Pick now any $0\leq u <\mu(k)$.
$\HF_{R/I} (k, u)$ is the number of ways one can choose the position of $k$ queens and $u$ squares unattacked by these queens.
For each such choice, one can extend  the set of $u$ unattacked squares to a maximal set of $v$ unattacked squares by the same $k$ queens.
To obtain $\Phi(k,u)$ we need to count only those choices for which $u=v$ or, equivalently,
we need to subtract from $\HF_{R/I} (k, u)$ the number of configurations which which extend to a maximal one with $v>u$ unattacked squares.
The induction hypothesis implies that there are exactly $\Phi(k,v)$ configurations with $k$ queens and a maximal set of $v$ unattacked squares,
and each one of these produces $\binom{v}{u}$ configurations with $k$ queens and $u$ unattacked squares
which can be extended to a maximal set of $v$ unattacked squares. Subtracting all these, we get the desired result.
\end{proof}

\bigskip
Table 1 lists the values of $\Phi(k,u)$ when $n=8$ for $3 \leq k \leq 43$ and $1 \leq u \leq 25$ (blank entries are zero.)
For example, the table shows that $\mu(8)=11$ and that $\Phi(8,\mu(8))=48$, which means that the largest number of unattacked squares one can have
when 8 queens are placed on a regular chessboard is 11, and that there are 48 such configurations.
This is the answer to a question
originally published by W.~W.~Rouse Ball in 1896 \cite{B} (see also chapter 34 in \cite{G}.)
This calculation was produced by {\tt FreeSquares}, a C++ program which can be found in \cite{K1}.
(There are several widely used computer packages which can compute multi-graded Hilbert series,
but unfortunately they are not very efficient.)

\begin{figure}[p]
\leavevmode
\includegraphics[height=10in]{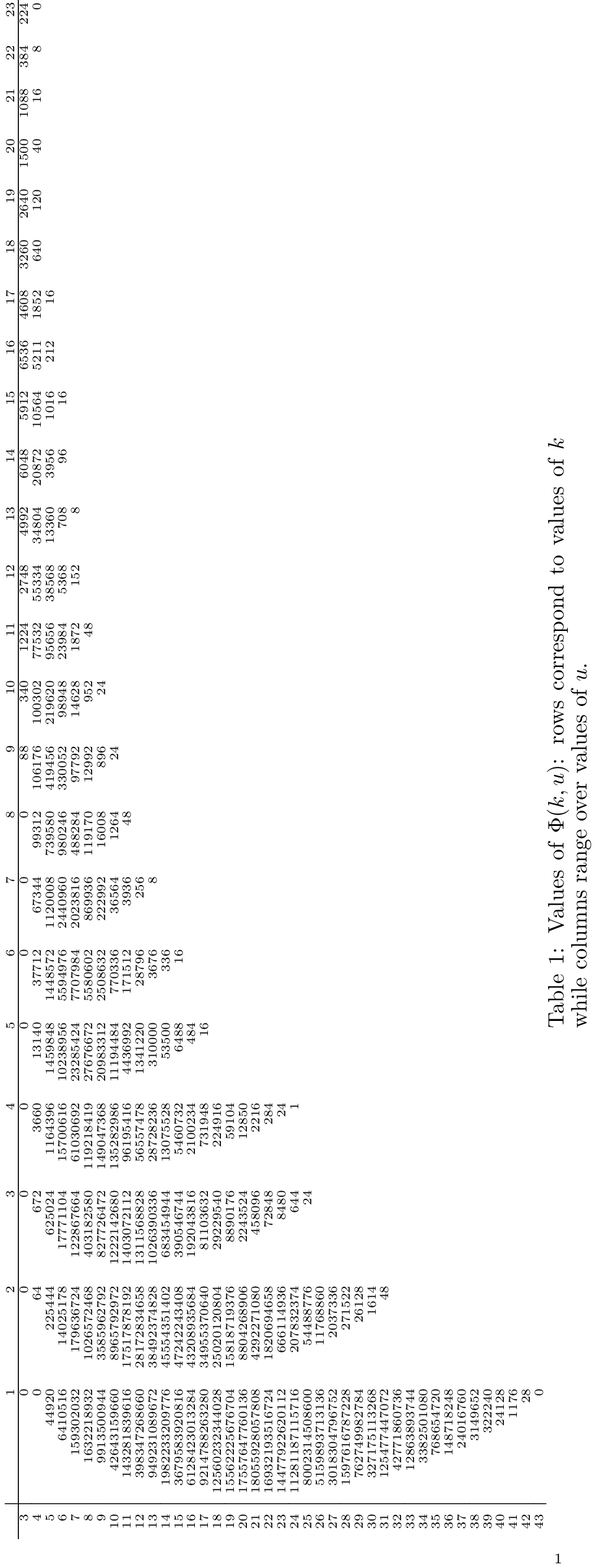}\\
\end{figure}

\bigskip
The method introduced in this section generalizes naturally to deal with graph-theoretical problems which we now describe.
Let $G$ be a finite graph. If $U$ and $W$ are disjoint sets of vertices of $G$ we say that
$U$ and $W$ are independent if there is no edge connecting a vertex in $W$ with a vertex in $U$.
For a given $k$ what is the maximal size of a set of vertices which is independent of a set of $k$ vertices?
In how many ways can one choose independent $U$ and $W$ with given size?

Let $\left\{ v_1, \dots, v_N \right\}$ be the vertices of $G$.
One obtains the solution to this more general problem by replacing
the ring $R$ with $K[x_1, \dots, x_N, y_1, \dots, y_N]$ and
the ideal $I$ above with the ideal generated by the squares of all the variables and
$$  \left\{ x_{i} y_{j} \,|\, (v_i,v_j)\mathrm{\ is\ an\ edge\ in\ } G \right\} .$$

\section{Knight moves in an infinite chessboard.}

We now consider the second set of questions mentioned in the introduction:
How many squares can a knight in an infinite chessboard reach in $d$ moves?
How many squares can be reached in $d$ moves and no less moves?
We will denote the first number with $f(d)$ and the second with $g(d)$.

The implementation of the results in this section relies on Gr\"obner bases techniques--
the reader may want to consult \cite{AL} for an introduction to Gr\"obner bases.
However, to appreciate the general ideas behind the approach of this section no knowledge of Gr\"obner bases is needed.

We again let $K$ be any field and let $R$ be the
$K$-subalgebra of $K[x_1,x_2,x_1^{-1},x_2^{-1}]$ generated
by
$$M=\left\{ x_1 x_2^2, x_1^2 x_2, x_1^{-1} x_2^2, x_1^{-2} x_2,
x_1 x_2^{-2}, x_1^2 x_2^{-1}, x_1^{-1} x_2^{-2}, x_1^{-2} x_2^{-1}\right\}.$$

The first step towards the solution of this problem is to realize that
$f(d)$ is the cardinality of
$M^d:=\{a_1 \dots a_d | a_1,\dots,a_d\in M\}$
while $g(d)$ is the number of elements
in $M^d$ but not in any $M^i$ for $i<d$.

We can produce a presentation for $R$ by mapping a polynomial ring
$S=K[y_1,\dots,y_8]$ to $R$ by $y_i \rightarrow m_i$ where $m_i$ is
the $i$th element of $M$.
We denote this mapping with $\Psi$.
Notice that the restriction of $\Psi$ to the set of
degree-$d$ monomials in $S$ gives a surjection onto the elements of $M^d$.

Let $\kappa$ be the kernel of the map above.
This kernel can be computed effectively using Gr\"obner bases techniques as follows:
let $I$ be the ideal of $k[u,x_1,x_2,y_1,\dots,y_8]$ generated by
\begin{gather*}
\{
u x_1 x_2-1, y_1-x_1 x_2^2, y_2-x_1^2 x_2, y_3 x_1-x_2^2, y_4 x_1^2-x_2,\\
y_5 x_2^2-x_1, y_6 x_2-x_1^2, y_7 x_1 x_2^2-1, y_8 x_1^2 x_2-1 \}
\end{gather*}
and fix an elimination order where $u,x_1,x_2$ are the largest variables.
Then $\kappa$ is generated by the elements of a Gr\"obner basis for $I$ which
do not contain the variables $u,x_1,x_2$ (cf. chapter 1 of \cite{Stu}.)
Recall also that $\kappa$ is a binomial ideal.

Notice that the ring $R$ is not very interesting: it is in fact identical to $K[x_1,x_1^{-1},x_2,x_2^{-1}]$
(here is a chess proof: $x_1\in R$ because a knight can move one square to the right in three moves. By symmetry also $x_1^{-1},x_2,x_2^{-1}\in R$.)
However, $S/\kappa$ is far more interesting for reasons explained below.

Since the restriction of $\Psi$ to the set of degree-$d$ monomials in $S$ is a surjection
onto $M^d$,
to find $f(d)$ we need to find the size of a maximal set of degree-$d$ monomials in $S$
which are distinct modulo $\kappa$.
Two such monomials $\mathbf{y}^\alpha$ and $\mathbf{y}^\beta$ are distinct modulo $\kappa$ if and only if
$\mathbf{y}^\alpha-\mathbf{y}^\beta$ is not in the largest homogeneous sub-ideal $H$ of $\kappa$.
It is easy to compute $H$:
the elements of $H$ are the elements of the homogenization of $\kappa$ with
respect to a new variable, say $t$, which do not involve $t$, thus we
can compute $H$ by homogenizing a Gr\"obner basis for $K$ using
a graded lexicographic order (cf. exercise 1.6.19 in \cite{AL}) and eliminating the variable $t$.
We notice that this Gr\"obner basis can be chosen to consist of binomials, and so $H$ is also a binomial ideal.

So we have reduced the problem of computing  $f(d)$
to the problem of finding the size of a maximal set of degree-$d$ monomials in $S$ which are distinct modulo $H$.
Fix any term ordering in $S$  and let $\mathcal{H}$ be a Gr\"obner basis for $H$ consisting of binomials.
Now for any two monomials $\mathbf{y}^\alpha > \mathbf{y}^\beta$ of the same degree,
$\mathbf{y}^\alpha\equiv \mathbf{y}^\beta$ modulo $H$ if and only if
$\mathbf{y}^\alpha$ reduces to $\mathbf{y}^\beta$ with respect to $\mathcal{H}$.
Since each reduction of a monomial with respect to $\mathcal{H}$ produces a new monomial (of same degree),
to produce a maximal set of degree-$d$ monomials in $S$ which are distinct modulo $H$
we may pick all monomials of degree $d$ which are non-zero modulo $\initial(H)$, i.e.,
$$f(d)=\dim_K \left(S/\initial(H)\right)_d=\dim_K \left(S/H\right)_d= \HF_{S/H} (d) $$
where the second equality is a celebrated theorem proved by F.~S.~Macaulay in \cite{M}.

An easy computation with Macaulay2 (\cite{GS})  shows that
$$\HS_{S/H}(t)=\frac{1+5t+12t^2-8t^4+4t^5}{(1-t)^3}$$
and that the Hilbert polynomial of $S/H$ is $1+4d+7d^2$.
Since
$$\HS_{S/H}(t) - \sum_{d=0}^\infty (1+4d+7d^2)t^d=-4t^2-4t$$
we obtain
$$f(d)=
\left\{
\begin{array}{l l}
1 & d=0 \\
8 & d=1 \\
33 & d=2 \\
1+4d+7d^2 & d\ge 3
\end{array}
\right.
$$

\bigskip
We now proceed to compute $g(d)$.
We again fix a monomial ordering in $S$ which refines the total degree ordering.
List all the monomials in $S$ in ascending order,
and let $B$ be the set of all degree-$d$ monomials in $S$ which are not congruent modulo $\kappa$ to a monomial appearing earlier in the list.
We now show that $g(d)=\#B$.

If for two distinct degree-$d$ monomials $\mathbf{y}^\alpha>\mathbf{y}^\beta$ we have
$\Psi(\mathbf{y}^\alpha)=\Psi(\mathbf{y}^\beta)$ then $\mathbf{y}^\alpha-\mathbf{y}^\beta\in \kappa$ contradicting the choice of $B$.
Hence the restriction of $\Psi$ to $B$ is injective.
Similarly, if for some degree-$d$ monomial $\mathbf{y}^\alpha$ there exist
a monomial $\mathbf{y}^\beta$ of degree $i<d$ so that $\Psi(\mathbf{y}^\alpha)=\Psi(\mathbf{y}^\beta)$ then $\mathbf{y}^\alpha-\mathbf{y}^\beta\in \kappa$ and since
$\mathbf{y}^\alpha>\mathbf{y}^\beta$ we get a contradiction to the choice of $B$.
Hence the restriction of $\Psi$ to $B$ is a surjection onto $M^d\setminus \cup_{i<d} M^i$.

Using the fact that $\kappa$ has a Gr\"obner basis generated by binomials we may deduce that
$B$ is the set of all monomials which are not in $\initial \kappa$ and so
$$g(d)=\dim_K \left(S/\initial(\kappa)\right)_d= \HF_{\left(S/\initial(\kappa)\right)} (d) .$$

Another straightforward computation with Macaulay2 shows that
$$\HS_{\left(S/\initial(\kappa)\right)}(t)=\frac{1+6t+17t^2+12t^3-8t^4-4t^5+4t^6}{(1-t)^2}$$
and that the Hilbert polynomial of $S/\initial(\kappa)$ is $28d-20$.
Since
$$\HS_{\left(S/\initial(\kappa)\right)}(t) - \sum_{d=0}^\infty (28d-20)t^d=4t^4+4t^3-4t^2+21$$
we obtain
$$g(d)=
\left\{
\begin{array}{l l}
1 & d=0 \\
8 & d=1 \\
32 & d=2 \\
68 & d=3 \\
96 & d=4 \\
28d-20 & d\ge 5
\end{array}
\right.
$$

The methods of this section also generalize in a natural way.
Let
$$W=\left\{
\left(\begin{array}{c} w_{1 1}\\ \vdots \\ w_{1 m} \end{array}\right),
\dots ,
\left(\begin{array}{c} w_{N 1}\\ \vdots \\ w_{N m} \end{array}\right)
\right\}
\subset \mathbb{Z}^m$$
be a finite set
and consider an infinite directed graph $G$ whose
vertex set is $\mathbb{Z}^m$ and for any $u,v \in \mathbb{Z}^m$, $\overrightarrow{(u,v)}$ is a directed edge
if and only if $v-u \in W$.

By replacing $R$ above and its presentation  $S\rightarrow R$ with the presentation
$$K[y_1, \dots, y_N] \rightarrow K\left[x_1^{w_{1 1}} \cdot \dots \cdot x_m^{w_{1 m}}, \dots , x_1^{w_{N 1}} \cdot \dots \cdot x_m^{w_{N m}}\right]$$
which maps $y_i$ to $x_1^{w_{i 1}} \cdot \dots \cdot x_m^{w_{i m}}$ for all $1\leq i\leq N$, we can, by following exactly the same procedures
as before, produce closed formulas for the functions
$f(d)$ which count how many endpoints all length $d$ paths starting at a fix vertex have,
and closed formulas for the functions $g(d)$ which count
how many vertices are at a distance of $d$ from a fixed vertex.

\begin{thm}
For any directed graph $G$ as above,
there exist polynomials $P(d)$ and $Q(d)$ so that $f(d)=P(d)$ and $g(d)=Q(d)$ for all $d \gg 0$.
\end{thm}
\begin{proof}
This is an immediate consequence of the fact that Hilbert functions are of polynomial type.
\end{proof}

\section*{Appendix: A Macaulay2 implementation.}

All the methods in this paper are easy to implement with existing computer systems.
As an example aimed to tempt the reader to experiment with these systems we present
a Macaulay2 program for the solution of the enumeration problem in the previous section:

\bigskip
\def\baselinestretch{0.8}
\begin{verbatim}
R=ZZ/101[u,a,b,y_{1}..y_{8},MonomialOrder=>Lex];
I={u*a*b-1_R,y_{1}-a*b^2,y_{2}-a^2*b,y_{3}*a-b^2,y_{4}*a^2-b,
y_{5}*b^2-a,y_{6}*b-a^2,y_{7}*a*b^2-1_R,y_{8}*a^2*b-1_R};
G=gens gb ideal I;
J=selectInSubring(3,G);

S1=ZZ/101[y_{1}..y_{8},t];
J=substitute(J,S1);
H0=homogenize(gens gb J,t);

S2=ZZ/101[t,y_{1}..y_{8},MonomialOrder=>Lex];
H0=substitute(H0,S2);
G=gens gb ideal H0;
H=selectInSubring(1,G);

S=ZZ/101[y_{1}..y_{8}];
J=substitute(J,S);
H=substitute(H,S);
print(hilbertSeries coker J);
print(hilbertPolynomial(coker J, Projective=>false));
print(hilbertSeries coker H);
print(hilbertPolynomial(coker H, Projective=>false));
\end{verbatim}

This produces the following output:

\begin{verbatim}
   6    5    4     3     2
4$T -4$T -8$T +12$T +17$T +6$T+1
--------------------------------
                   2
            (-$T+1)
28$i-20
   5    4     2
4$T -8$T +12$T +5$T+1
---------------------
              3
       (-$T+1)
   2
7$i +4$i+1
\end{verbatim}
\def\baselinestretch{1}


\end{document}